\documentclass{article}
\usepackage{amsmath}
\usepackage{amsthm}
\usepackage{epsf}

\def\N{\mathbf{N}}
\def\Z{\mathbf{Z}}
\def\QED{\unskip~\vrule width 5pt\smallskip}

\let\after=\circ

\newcommand\inv{^{-1}}
\newcommand\eps{\varepsilon}

\newtheorem{theorem}{Theorem}
\newtheorem{corol}[theorem]{Corollary}
\newtheorem{prop}[theorem]{Proposition}

\newcommand\vcsanepsf[2]{\hbox{$\vcenter{\hbox{#1\epsfbox{#2}}}$}}

\title{Some simple bijections involving lattice walks and ballot sequences}

\author{Marc A. A. van Leeuwen%
        \thanks{Partially supported by NSF grant DMS-0554278}
        \\
        Universit\'e de Poitiers}

\date{}

\begin{document}

\maketitle

\begin{abstract}
In this note we observe that a bijection related to Littelmann's root
operators (for type~$A_1$) transparently explains the well known enumeration
by length of walks on~$\N$ (left factors of Dyck paths), as well as some other
enumerative coincidences. We indicate a relation with bijective solutions of
Bertrand's ballot problem: those can be mechanically transformed into
bijective proofs of the mentioned enumeration formula.
\end{abstract}

\section{Introduction}

When considering a combinatorial formula that can be interpreted as equating
the outcomes of two (families of) enumeration problems, a proof in the form of
a bijection (or family thereof) between the sets in question is often
considered to be of more value than one based on other methods, such as
manipulations of formal power series. Whether this is really the case depends
to a large extent on the nature of the bijection. The best situation is one
where the bijection can be interpreted as simply relating two ways describing
the same underlying combinatorial object; for instance the bijection between
Dyck paths and balanced sequences of parentheses, which simply interprets
up-steps as `{\tt(}' and down-steps as `{\tt)}', is of this nature. A
bijection that consists of transforming an object of the first kind into an
object of the second kind by a single traversal performing some kind of
substitution, and in such a way that the relation between input and final
output can be easily perceived, will do almost just as well. In both cases it
is often possible to relate one or more additional statistics on the input
with an statistics on the output, and thus refine a simple identity of natural
numbers into an one of polynomials with non-negative coefficients (the idea of
$q$-analogues is largely based on this principle). However when the bijection
is based on a more complicated algorithm, or is obtained by a composition of
several bijections, then the situation may become much less transparent, to
the point of providing no more insight (possibly even less) than a proof by
formal algebraic manipulations.

The dividing line between simple and complicated bijections is not always
clear, and a bijection may be rendered more transparent by a particular point
of view. For instance, one can map binary plane trees with $n$ internal nodes
(and $n+1$ leaves) to Dyck paths of length $2n$ by traversing the tree in
pre-order, recording an up-step for each internal node encountered and a
down-step for every leaf except the final one (cf.~\cite{Stanley EC2},
proposition 6.2.1(i),(ii)), and such a path~$P$ can be further transformed
into a plane tree with $n+1$ vertices, namely one whose root node has
descendents that correspond (in order) to the minimal Dyck-path factors
of~$P$, and are obtained by recursively transforming those factors after
removing their two extremal steps. The resulting bijection from binary plane
trees to plane trees appears to be rather opaque, and not provide much more
insight than the observation that the generating series for both enumeration
problems satisfy (for simple reasons) the quadratic relation $C=1+XC^2$ that
characterises the series $C=\sum_nC_nX^n$ of the Catalan numbers, which
numbers therefore solve both problems. However, from the point of view of a
programming language like LISP (in which non-empty lists are represented by a
binary node with links pointing to the first element and to the remainder of
the list), one can interpret the binary tree as the internal representation of
the recursively nested list of lists corresponding to the planar tree (of
which each node is interpreted as the list of its descendents). This point of
view makes the correspondence, at least to us, much more transparent.

On the other hand, when a bijection is defined by an algorithm that involves
the repetition of some operation a (finite but) variable number of times,
until reaching some desired condition, then this bijection is likely to be
quite opaque. This is for instance probably the reason that use of the
``involution principle'' (due to Garsia and Milne) is generally considered to
be less desirable in bijective proofs. An example where such iteration can
arise, is when one is given a bijection $f:A\to B$ where $A,B$ are finite
subsets of some set~$X$, and one deduces from it a bijection
$g:X\setminus{B}\to X\setminus{A}$ between their complements, by iterating $f$
as often as possible; in other words $g(x)=f^n(x)$ for the smallest~$n\in\N$
for which $f^n(x)\notin{A}$ (which must exist because all other $f^i(x)$ are
distinct elements of~$A$). Whether the bijection $g$ obtained by this
``complementation principle'' is transparent at all, depends on whether the
effect of iterating~$f$ can be easily understood, and notably on whether it is
easy to predict the number~$n$ of iterations that can be applied to a given
element~$x$.
\unskip\footnote {In fact the involution principle can be seen as an
  application of the complementation principle to the situation of a ``signed
  set''~$Y$ consisting of disjoint union of a set~$X$ of ``positive'' elements
  and a finite set~$N$ of ``negative'' elements, equipped with injections
  $f_1,f_2:{N\to X}$; one can obtain a bijection $X\setminus f_1(N)\to
  X\setminus f_2(N)$ by taking $A=f_2(N)$, $B=f_1(N)$, and $f=f_1\after
  f_2\inv$. The injections can be extended to fixed-point free involutions of
  the union of~$N$ and its image, and then to involutions of~$Y$ by fixing the
  remaining points of~$X$; this explains the name of the principle. In
  practice the fixed-point free involutions in this description are usually
  directly obtained as sign-reversing partial involutions of~$Y$ that are
  always defined at negative elements, and which are then artificially
  extended as we did to involutions of~$Y$.}

In this note we will consider the particular case of some well known
enumerative results involving lattice paths, or equivalently walks on the
one-dimensional lattice~$\Z$. Our walks will always start at~$0$, and we first
consider the ``binary'' case where each step changes the position either
by~$+1$ or by~$-1$. Each walk can be transformed into a lattice path in (the
``diagonal'' index~$2$ sub-lattice of)~$\Z^2$ starting at the origin, in which
steps advancing by~$\eps$ become path segments advancing by $(1,\eps)$ (which
we draw with the first ``time'' coordinate increasing downwards). We shall
switch between these points of view whenever convenient.

We call a finite walk
``recurrent'' if it ends at~$0$, and ``positive'' if it is a walk on~$\N$
(``non-negative'' would be more precise, but tiresome). Viewed as lattice
paths, (binary) recurrent positive walks correspond precisely to Dyck paths.

The most basic case of the enumerative coincidences that we shall study is the
fact that there are $\binom{2n}n$ positive walks of length~$2n$, a number that
also (and more obviously) counts the recurrent walks of that length. This
result appears to be well known, at least in the lattice path community, but
in view of its simplicity it is somewhat surprising that it does not receive
prominent mention in the enumerative combinatorics literature. We do not know
whether any nice bijective proofs for this result are known, but it would at
least seem that none are ``well known''. This note proposes a simple bijective
proof that, although it involves the iteration of an operation a varying
number of times to transform recurrent walks into positive ones or vice versa,
is about as transparent as one could wish for; notably the number of
iterations required can be immediately read off from the initial walk, and it
is possible to transform the initial to the final path in a single ``pass''
along the path. This bijection also allows giving a bijective proof of the
generating series identity $(\sum_{n\in\N}\binom{2n}nX^n)^2=\frac1{1-4X}$, for
which there does not appear to be an equally transparent proof using either of
the above interpretations of the number~$\binom{2n}n$ individually.

The result can be slightly generalised, while essentially keeping the same
proof, in a few ways. One can drop the restriction to walks of even length
provided the ``recurrent'' requirement is relaxed to ending either at~$0$ or
at~$1$ (since clearly being recurrent is a tall order for odd-length walks).
One can also consider ``ternary'' walks by allowing steps that stay in place
(so that the paths corresponding to recurrent positive walks are Motzkin
paths), in which case the positive walks are still in bijection with the
``almost recurrent'' walks, those that end either at~$0$ or at~$1$ (the nature
of our proof will make clear why one must leave a unit of freedom for the
ending point of the walk). Finally one can formulate a corresponding result
for $n$-dimensional walks with a fairly large choice for the set of basic
steps allowed: no coordinate should be allowed to change by more than a unit
at a time, and the set should be symmetric with respect to negation of each of
the coordinates individually.

The bijection we propose is defined by iterating a basic ``raising'' operation
as often as possible. We do not consider it as obtained from the
complementation principle mentioned above (although it can be), but rather
from a telescoping sum of identities, each of which is closely related to the
famous ``ballot problem'' of J.~Bertrand, \cite{Bertrand}. Our bijective proof
gives rise to a bijective solution of that problem, which appears to be new;
at least it is different both from the original proof of
D.~Andr\'e~\cite{Andre} (which performs a cyclic rearrangement of votes), and
from proofs based on diagonal reflection of part of a lattice path (the
``reflection method''). However we shall see that conversely any bijective
proof of the ballot problem can be iterated so as to obtain a proof of the
identities relating positive and recurrent walks, of which we therefore obtain
several different ones. Comparing these, we find that some proof methods that
transparently solve the ballot problem lead to rather opaque bijective proofs
after iteration.

\section{Positive and (almost) recurrent walks}

A very basic kind of lattice walks is that of walks on the one dimensional
lattice~$\Z$, starting at~$0$ and moving a unit in either direction at each
step. The parity of the point reached after $n$ steps is necessarily that
of~$n$, and the statistic of the end point on the set of the $2^n$ such walks
gives a binomial distribution on the points of the required parity. Another,
more restricted, class of walks that we shall consider is that of walks on the
subset~$\N$ of the one dimensional lattice, still starting at~$0$ and moving a
unit in either direction at each step; these form the subset of the walks
on~$\Z$ that never visit the value~$-1$. For comparison, here is an initial
portion of Pascal's triangle, displayed in the usual fashion with rows
symmetrically growing as one moves downwards,
\unskip\footnote {This is not the way Pascal drew his ``Triangle
  Arithmetique''; he used a horizontal grid of cells and the
  rule ``Le nombre de chaque cellule est egal {\`a} celuy de la cellule qui la
  precede dans son rang perpendiculaire, plus {\`a} celuy de la cellule qui la
  precede dans son rang parallele.''~\cite{Pascal}}
and a corresponding array of numbers counting walks on~$\N$.
\[\setbox0=\hbox{00}\dimen0=\wd0
  \vtop{\halign {&\hbox to\dimen0{\hss#\hss}\cr
   &   &   &  &  &  &  &  &  &  & 1\cr
   &   &   &  &  &  &  &  &  & 1&  & 1\cr
   &   &   &  &  &  &  &  & 1&  & 2&  & 1\cr
   &   &   &  &  &  &  & 1&  & 3&  & 3&  & 1\cr
   &   &   &  &  &  & 1&  & 4&  & 6&  & 4&  & 1\cr
   &   &   &  &  & 1&  & 5&  &10&  &10&  & 5&  & 1\cr
   &   &   &  & 1&  & 6&  &15&  &20&  &15&  & 6&  & 1\cr
   &   &   & 1&  & 7&  &21&  &35&  &35&  &21&  & 7&  & 1\cr
   &   &  1&  & 8&  &28&  &56&  &70&  &56&  &28&  & 8&  & 1\cr
   &  1&   & 9&  &36&  &84& &126& &126&  &84&  &36&  & 9&  & 1\cr
  1&   & 10&  &45& &120& &210& &252& &210& &120&  &45&  &10&  & 1\cr
}}
\quad
  \vtop{\halign {&\hbox to\dimen0{\hss#\hss}\cr
  & 1\cr
 0&  & 1\cr
  & 1&  & 1\cr
 0&  & 2&  & 1\cr
  & 2&  & 3&  & 1\cr
 0&  & 5&  & 4&  & 1\cr
  & 5&  & 9&  & 5&  & 1\cr
 0&  &14&  &14&  & 6&  & 1\cr
  &14&  &28&  &20&  & 7&  & 1\cr
 0&  &42&  &48&  &27&  & 8&  & 1\cr
  &42&  &90&  &75&  &35&  & 9&  & 1\cr
}}
\]
Note that the second array obeys the same recurrence relation that the first
one does (every ``internal'' entry is the sum of those directly above it),
only the boundary condition is changed, namely by requiring the entries in the
column of~$-1$ (the leftmost ones displayed) to be~$0$, reflecting the fact
that walks that would visit negative numbers are excluded. In fact the second
array can be obtained from the first by subtracting from it a copy of itself
that is shifted two units (the distance between adjacent entries) to the left,
which produces values~$0$ in the column~$-1$ for symmetry reasons; one retains
the part to the right of that column. It follows that the sum of the entries
in any row of the second array is equal to the (most) central entry in the
corresponding row of the first array; in the terminology of the introduction,
there are as many walks of a given length that are recurrent or (in the case
of odd length) almost recurrent (ending at~$1$) as there are positive walks.
It is this ``coincidence'' that we wish to bijectively explain in this note.
We can illustrate the classes of paths between which we seek a bijection
graphically as follows (the drawn paths are just examples):
\[
  \epsfysize=.4\hsize % this one is actually square; the x-size is relevant
  \epsfbox{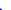} \qquad
  \epsfysize=.4\hsize % but this one is not square; y-sizes must match though
  \epsfbox{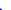}
\]
The way in which we obtained the equality of the number of paths of these two
types is rather simple, and one may seek to obtain a bijection from it by
applying general principles; this will be discussed in the next section. Here
however we shall directly consider a correspondence between all positive walks
and recurrent walks, without first constructing one corresponding to the fact,
implicitly mentioned in the reasoning above, that for any
$\frac{n}2\leq{k}\leq{n}$ there are $\binom{n}k-\binom{n}{k+1}$ positive walks
that end at the value~$2k-n$.

If a recurrent walk~$w$ happens to be positive as well (it corresponds to a
Dyck path), then it can be made to correspond to itself, without applying any
operations. In the contrary case there will certainly be a first down-step
in~$w$ that reaches~$-1$, and then later on possibly a first down-step that
reaches~$-2$, and so forth. All these down-steps that for the first time reach
a given negative number are changed into up-steps to form the positive
walk~$w'$ corresponding to~$w$; if there were $d$ such steps (so that $-d$ is
the most negative number that $w$ visits), then~$w'$ will end at the
number~$2d$. Any positive walk of even length ends in an even non-negative
number; the unique recurrent walk~$w$ corresponding to it can be found by
setting~$d$ to half that final number, and changing the $d$ up-steps that
immediately follow the last visit to respectively the numbers $0$,~$1$,
\dots,~$d-1$ into down-steps. This correspondence can be extended
straightforwardly to include odd-length walks, making those that end at~$1$
and reach~$-d$ as most negative number correspond to positive walks ending
at~$2d+1$.

Note that the bijection establishes a fact that was not evident in our
original argument, namely that $\binom{n}k-\binom{n}{k+1}$ not only counts the
positive walks that end at~$2k-n$, but also the walks of ``depth''
$k-\lceil\frac{n}2\rceil$ ending at $0$ or~$1$.

Before we state more formally the result thus obtained, we shall generalise it
slightly by allowing in addition to up-steps and down-steps also neutral
steps, which stay at the same point. Both our initial reasoning and the
construction of a correspondence remain valid without much modification for
these more general walks, although of course the numbers of walks increase,
and there is no longer a parity condition for the end point of the walks.
Instead of Pascal's triangle and it anti-symmetrised counterpart we obtain
as arrays of numbers the coefficients of $(X\inv+1+X)^n$, somewhat ambiguously
called trinomial coefficients:
\[\setbox0=\hbox{000,}\dimen0=\wd0
  \vtop{\halign {&\hbox to\dimen0{\hss#\hss}\cr
    &   &   &   &   &   &   &  1\cr
    &   &   &   &   &   &  1&  1&  1\cr
    &   &   &   &   &  1&  2&  3&  2&  1\cr
    &   &   &   &  1&  3&  6&  7&  6&  3&  1\cr
    &   &   &  1&  4& 10& 16& 19& 16& 10&  4&  1\cr
    &   &  1&  5& 15& 30& 45& 51& 45& 30& 15&  5& 1\cr
    &  1&  6& 21& 50& 90&126&141&126& 90& 50& 21& 6& 1\cr
   1&  7& 28& 77&161&266&357&393&357&266&161& 77&28& 7& 1\cr
}}
\]
and
\[\setbox0=\hbox{000,}\dimen0=\wd0
  \vbox{\halign {&\hbox to\dimen0{\hss#\hss}\cr
&&&&&&  &  1\cr
&&&&&& 0&  1&  1\cr
&&&&&& 0&  2&  2&  1\cr
&&&&&& 0&  4&  5&  3&  1\cr
&&&&&& 0&  9& 12&  9&  4& 1\cr
&&&&&& 0& 21& 30& 25& 14& 5& 1\cr
&&&&&& 0& 51& 76& 69& 44&20& 6& 1\cr
&&&&&& 0&127&196&189&133&70&27& 7& 1\cr
}}.
\]
The entries of the second array are differences of trinomial coefficient two
places apart, so each of its row sums is given by the corresponding middle
trinomial coefficient \emph{plus} one its neighbours. Like for the first
correspondence, we illustrate graphically the types of paths matched by the
bijection, and as in the first case the paths depicted actually match under
our correspondence.
\[
  \epsfysize=.4\hsize % this one is almost square; the x-size is relevant
  \epsfbox{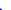} \qquad
  \epsfysize=.4\hsize % but this one is not square; y-sizes must match though
  \epsfbox{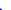}
\]

We can now formulate our main result.

\begin{theorem}
Within the class of walks on~$\Z$ starting at~$0$ and with steps advancing by
$+1$, $0$ or $-1$, there is a bijection, conserving both the length of the
walk and the number of steps~$0$, between on one hand the walks that end
either at~$0$ or at~$1$, and on the other hand the walks that do not visit
negative numbers. The bijection maps walks ending at~$e\in\{0,1\}$ and whose
minimal number visited is~$-d$, to walks ending at~$2d+e$, and is realised by
reversing the direction of the $d$~down-steps that first reach respectively the
numbers $-1$, $-2$, \dots, $-d$.
\end{theorem}

\emph{Proof.} Calling ``Motzkin walk'' any sub-walk starting and ending at a
same number~$m$ while not visiting any number less than~$m$ (these correspond
to Motzkin sub-paths in the lattice path point of view), any walk in the class
considered ending at~$e\in\Z$ and whose minimal number visited is~$-d$ can be
uniquely written as a composition of $2(2d+e)+1$ sub-walks that are
alternatingly a (possibly empty) Motzkin walk and a single non-stationary
step, those $2d+e$ single steps being (in order) $d$ down-steps and $d+e$
up-steps; they are precisely the steps not contained in any Motzkin walk. The
domain and codomain of the bijection are characterised by $e\in\{0,1\}$
respectively by $d=0$, and the bijection, which reverses the $d$ down-steps in
this decomposition, produces the same Motzkin walk factors, changing the
parameters from $(d,e)$ to~$(0,2d+e)$. \QED

Although we have described the bijection as a single transformation, it can be
obtained by repeating a same operation, which reverses only the direction of a
single step, $d$~times in succession. Doing so, there is no choice but to
start reversing the last one of those $d$ steps, the one that first attains
the global minimum~$-d$ of the walk: reversing any of the other ones would
result in a step that becomes part of a Motzkin walk without in general any
means to tell from that walk alone which step it was. The final step however
remains outside any Motzkin walks after reversal, and in fact becomes the last
up-step that starts at the global minimum of the modified walk, which has
become~$-d+1$; this description shows that the original walk can be
reconstructed given only the modified walk. Now repeating the operation of
reversing the down-step that first attains the global minimum of the walk will
successively reverse the $d$ steps indicated in our bijection, in reverse
order of appearance in the walk; after this the operation cannot be further
repeated because the global minimum has become~$0$, and no down-step is needed
to first attain it. The successive steps of the transformation are illustrated
in figure~1, in a case with $d=3$.

\begin{figure}[hb]
\def\map{\hfil\longrightarrow\hfil}
\hbox to \hsize
{\hfil\hfil\hfil\hfil$
  \vcsanepsf{\epsfysize=52mm}{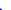} \map
  \vcsanepsf{\epsfysize=52mm}{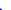} \map
  \vcsanepsf{\epsfysize=52mm}{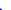} \map
  \vcsanepsf{\epsfysize=52mm}{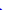} $\hfil\hfil\hfil\hfil}
\caption{From a recurrent walk to a positive walk by pushing minimum upwards}
\end{figure}

The reverse procedure then consists of repeating the following operation,
which decreases the end point of the walk by~$2$, until that end point is $0$
or~$1$: reverse the first up-step that starts at the global minimum of the
walk. Unlike the forward operation, this backward operation could in general
be repeated in a reasonable way even after the terminating condition is
reached. Doing so until no up-step of the given type can be found would result
in a bijection from the set of all positive walks to those ending at their
global minimum, and mapping those of the former kind ending at~$k$ to those of
the latter kind ending at~$-k$. There is more obvious bijection with the same
property, namely the one consisting of taking the steps of the walk in reverse
order and in the opposite direction (in terms of the corresponding lattice
paths, this corresponds to reflection in a horizontal line, and then shifting
horizontally to match the starting point).

The operation described on walks or paths is certainly not new. It appears
that its most prominent occurrence is in the representation theory of Lie
groups: it is an instance of Littelmann's root operator~$e_\alpha$ on paths
of~\cite{Littelmann}, for the most basic case of type~$A_1$. However, in
diverse settings and under various equivalent descriptions, this operation had
been known long before; notably it occurs in~\cite{Littlewood Richardson}
where it is used to prove the simplest case of the Littlewood-Richardson rule.

In those applications there can in fact be more than one such operation,
acting on paths in a space of higher dimension, and involving reflections
(applied to path segments) in different directions, namely the simple
reflections for a root system. It would be interesting to find enumerative
consequences of the operations in those settings, but it appears that there
are none that are easy to state if the simple reflections, and therefore the
associated root operators, do not commute. In fact, although iterating the
``raising operators'' does give, in spite of their non-commutation, a well
defined map from arbitrary paths to dominant ones (those that remain on the
positive side of each reflection hyperplanes considered), this map restricted
to for instance the recurrent paths is not injective (the size of the fibre
above a given dominant path ending at~$\lambda$ is the dimension of the
zero-weight space of a representation associated to~$\lambda$).

For this reason we shall state only a generalisation corresponding to the
type~$A_1^n$, where the reflections commute. In this case we can in fact
ignore the Lie theoretic point of view, which adds nothing that is not obvious
combinatorially. We can treat each of the~$n$ coordinate directions
separately, repeatedly applying reflections in that direction to certain
segments of a path until that coordinate is non-negative throughout the path.
These operations for different coordinates commute and can therefore be
applied independently. The set of basic steps from which the paths may be
constructed may however involve some dependence between the coordinate
directions (it could for instance insist that exactly two coordinates change
at each step), as long as it is symmetric with respect to each of the $n$
coordinate reflections, and no step involves a change outside the
set~$\{-1,0,1\}$ to any one coordinate. Therefore we state the result as
follows.

\begin{theorem}
Let any subset $S$ of $\{-1,0,1\}^n$ be fixed that is stable under each
reflection that negates a single coordinate. Then within the set of walks
on~$\Z^n$ starting at the origin and with steps in~$S$, there is a length
preserving bijection from those walks that end in a point of~$\{0,1\}^n$ (a
vertex of the unit hypercube) to those walks that remain at all times inside
the (weakly) positive orthant. It is defined by applying the bijection of
theorem~1 to each coordinate separately.
\end{theorem}

\section{Convolution of central binomial coefficients}

There is another enumerative identity for which the bijection we studied
provides insight. It is well known that the generating series
$S=\sum_{n\in\N}\binom{2n}nX^n$ of central binomial coefficients is given by
$S=\smash{\frac1{\sqrt{1-4X}}}$; this can for instance be established by
developing $\smash{(1-4X)^{-\frac12}}$ using the power series binomial formula
with arbitrary exponent (cf.~\cite{Stanley EC1}, Chapter~1, exercise~4).
However, the formula means that $S^2$ is a geometric series with ratio~$4X$
and constant term~$1$, so that one has
\begin{equation}
  \label{twopathid}
  \sum_{i+j=n}\binom{2i}i\binom{2j}j=2^{2n} \quad\mbox{for all $n\in\N$.}
\end{equation}
The question of explaining this identity combinatorially is an old one;
according to \cite{Stanley EC1} p.~52 and \cite{Sved} it was raised by
P.~Veress and solved by G.~Hajos in the 1930s.
% (a more precise reference would be welcome).
In spite of the simple form of (\ref{twopathid}), it is remarkably difficult
to do this if one interprets the summand in the most obvious way as counting
pairs of recurrent walks of lengths $2i$ and $2j$ respectively. Although many
answers have been proposed, often involving operations closely related to
those we have been discussing (see \cite{Sved}), none appear so transparent as
to really explain the simplicity of the identity. With what we have seen
above, we may however also interpret the central binomial coefficients as
counting the positive walks of the indicated lengths. Doing so for both
coefficients does not make the problem much easier, but if one interprets the
first one as counting recurrent walks and the second on as counting positive
walks, then a solution presents itself naturally.

From the concatenation of a recurrent walk and a positive walk, the factors
cannot in general be uniquely reconstructed, but this decomposition becomes
unique if one inserts a single up-step between the two factors, as this will
be the last step to start at~$0$. The result is a walk of odd length~$2n+1$
ending at a number~$>0$, and conversely any such walk has a well defined last
up-step starting at~$0$, allowing a decomposition of the indicated kind. Thus
we are led to interpret the second member $2^{2n}$ not as counting the set of
all walks of length~$2n$, but as counting the set of walks of length~$2n+1$
ending at a positive number. This is the most subtle twist: the latter set
clearly has $2^{2n+1}/2=2^{2n}$ elements as well, but is not in obvious
bijection with the former set.
\unskip\footnote {Bijections between these sets can be found, but do not add
  much understanding to the obvious fact that their numbers agree. One such
  bijection uses by cut-and-paste involving two unequal size parts:
  interpreting walks as ballot sequences as in the next section, one can map a
  sequence of length~$2n+1$ in which $A$ beats~$B$ to an arbitrary one of
  length~$2n$, by singling out and removing one particular vote (the first or
  the last one are obvious choices); if it is a vote for~$A$ then the
  remaining sequence (weakly favourable for~$A$) is returned, and otherwise
  its complement (all votes inverted) is returned, which is strictly
  favourable for~$B$.}
Thus in our opinion the key to
understanding (\ref{twopathid}) does not lie in finding a clever bijection,
but in choosing the correct interpretation of its expressions; the
bijection (concatenation with an up-step interposed) then becomes a
triviality.

Since the correspondence between recurrent walks and positive ones continues
to hold in the presence of neutral steps, we can generalise the
identity~(\ref{twopathid}) to this setting. We shall allow any fixed
number~$t$ of distinguished kinds of neutral steps, so that we shall recover
the previous binary case for $t=0$, and the ternary case (Motzkin type walks)
for $t=1$ (or one may consider that we are doing just the latter case, but
keeping track of the number of neutral steps in the exponent of~$t$, viewed as
an indeterminate). The number $R_{t,l}$ of all such walks that are recurrent
and of length~$l$ is then equal, due to theorem~1, to the number of positive
walks of the same length ending at an even number; the parity condition is due
to the fact that the basic raising operation advances the end point by~$2$,
and the positive walks ending at an odd number similarly correspond to the
almost recurrent walks ending at~$1$.

Given a pair consisting of a recurrent walk of length~$i$ and a positive walk
of length~$l-i$ ending at a positive number, one can concatenate them with a
single up-step inserted between them, to obtain a walk of length~$l+1$ ending
at a positive odd number. Again one obtains, by combining all cases $0\leq
i\leq l$, a bijection to the set of all walks of the latter type, as these can
be uniquely decomposed at their last (up-)step starting at~$0$. Counting all
walks of length~$l+1$ ending at a positive odd number is only slightly more
difficult than in the case $t=0$. Their number is half that of all walks of
that length ending at an odd number, while the number of all walks of that
length, and the difference between the number of those ending at even and odd
numbers, are respectively given by the evaluations at $X=1$ and $X=-1$ of
$(X\inv+t+X)^{l+1}$. Combining this one gets
\begin{equation}
  \sum_{i+j=l}R_{t,i}R_{t,j} = \frac{(t+2)^{l+1}-(t-2)^{l+1}}4
\end{equation}
where $R_{t,i}$ is the coefficient of~$X^0$ in~$(X\inv+t+X)^i$. One has
$R_{0,i}=0$ whenever $i$ is odd; this makes the equation trivial when $t=0$
and $l$ is odd, while for $t=0$ and $l=2n$ one recovers
equation~(\ref{twopathid}). For arbitrary~$t$ one can convert this equation
into one of generating series, summing the individual geometric series
$\sum_i(t\pm2)^{i+1}X^i$ obtained from the right hand side to
$\frac{t\pm2-(t^2-4)X^2}{1-2tX+(t^2-4)X^2}$, and simplifying to
\begin{equation}
\label{diagonaleq}
\Bigl(\sum_{i\in\N}R_{t,i}X^i\Bigr)^2 = \dfrac1{1-2tX+(t^2-4)X^2}.
\end{equation}
This leads in particular to the expression
$\sum_{i\in\N}R_{1,i}X^i=1/\sqrt{1-2X-3X^2}$ for the generating series of the
middle trinomial coefficients. This also explains why that ``diagonal''
series, and its relative~$S$ for the middle binomial coefficients, not only
satisfy a quadratic relation over the rational functions in~$X$
(cf.~\cite{Stanley EC2}, theorem~6.3.3), but are actually are square roots of
a rational function (cf.~[loc.~cit.] exercise~6.42), unlike the series for
Motzkin paths and Dyck paths, although the latter satisfy simpler recurrences.

The form of equation~(\ref{diagonaleq}) begs consideration of the case $t=2$
as well, for which we obtain $\sum_{i\in\N}R_{2,i}X^i=1/\sqrt{1-4X}$, the same
generating series as~$S$ above. This leads to the question: why is the number
of recurrent walks of length~$n$ with steps advancing $+1$ or $-1$ and two
distinguished kinds of neutral steps the same as the number of recurrent walks
of length~$2n$ with only steps advancing $+1$ or $-1$? We leave this as an
(easy) exercise to the reader.

\section{Relation with Bertrand's ballot problem}

Finally we want to make explicit a relation of our basic bijections with the
famous ``ballot problem'' of Joseph Bertrand~\cite{Bertrand}. The problem he
presented, and immediately solved, is to compute the conditional probability,
given that an election between persons $A$ and~$B$ is won by~$A$ with $m$ out
of~$\mu$ votes (so $2m>\mu$), that during the sequential counting of the votes
$A$ has had a strict lead over~$B$, from the counting of the first vote
onwards. The answer given, $\frac{2m-\mu}\mu$, is justified by an argument
counting the number $P_{m,\mu}$ of ``favourable'' ballot sequences as a
function of $\mu$ and $m>\frac\mu2$, stating that it satisfies the same
recurrence as the binomial coefficients~$\binom\mu{m}$ (but without mentioning
them), from which recurrence (and the implicit condition $P_{m,2m}=0$) the
general formula for~$P_{m,\mu}$ can be deduced. Ballot sequences correspond
straightforwardly to walks on~$\Z$, and the favourable ones are those that
start with an up-step (a vote for~$A$), and then never return to~$0$. Their
number equals that of the positive walks of length $\mu-1$ with $m-1$
up-steps; as we have seen this gives
$P_{m,\mu}=\binom{\mu-1}{m-1}-\binom{\mu-1}m$, and one must assume this is the
general formula Bertrand hinted at. Indeed after division by~$\binom\mu{m}$
one obtains $\frac m\mu-\frac{\mu-m}\mu=\frac{2m-\mu}\mu$ as claimed.

Apparently Joseph's solution method differs little from the one by which we
initially counted positive walks. Yet he remarks that ``it seems probable that
so simple a result could be proved by a more direct method''. This challenge
is taken up by D{\'e}sir{\'e} Andr{\'e}, who in the very same issue of the
\textsl{Comptes Rendues} proposes a solution based on a combinatorial
argument~\cite{Andre}. He renames the parameters as $\alpha=m$ and
$\beta=\mu-m$ (the respective numbers of votes for $A$ and~$B$), and bases his
argument on the observation that the proposed probability
$\frac{\alpha-\beta}{\alpha+\beta}$ means that the complementary probability
$\frac{2\beta}{\alpha+\beta}$ (for the cases where $A$ did not maintain a
strict lead) is due for exactly half of it to the possibility that $A$ already
failed to win the very first vote, and for the other half to the possibility
that $A$ obtains an initial lead, but fails to maintain it. To simplify our
discussion we shall call these parts of the unfavourable scenario respectively
``bad'' and ``ugly'', and the favourable scenario ``good''. The conditional
probability of the bad case, where $B$ wins the first vote, is clearly equal
to the proportion $\smash{\frac\beta{\alpha+\beta}}$ of votes cast for~$B$; in
fact Andr{\'e} just computes this probability as the quotient
$\binom{\alpha+\beta-1}\alpha/\binom{\alpha+\beta}\alpha$. The heart of his
proof is a bijection establishing equality of the numbers of bad and ugly
cases.

Representing ballot sequences, or walks on~$\Z$ by lattice paths (drawn as
before with path segments in downward diagonal directions), the three
scenarios considered are illustrated in figure~2.
\begin{figure}[htb]
\hbox to \hsize
{\vcsanepsf{\epsfysize=37mm}{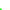} \hss
 \vcsanepsf{\epsfysize=37mm}{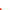} \hss
 \vcsanepsf{\epsfysize=37mm}{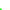}}
\caption{The good, the bad and the ugly cases}
\end{figure}
The equality observed by Andr{\'e} means the following: among the lattice
paths across a given rectangle, there are as many that start with a step along
the short side as there are that start with a step along the long side and
then later return at least once to the line bisecting the angle of the
rectangle at the starting point. The number of good cases, corresponding to
lattice paths from $(1,1)$ (which we take to be on the long side) to
$(\alpha+\beta,\alpha-\beta)$ that avoid the indicated line, can therefore be
computed as the total number of cases minus twice the number of bad cases,
corresponding to lattice paths from $(1,-1)$ to $(\alpha+\beta,\alpha-\beta)$;
this gives
\begin{equation}
  \binom {\alpha+\beta}\alpha - 2\binom{\alpha+\beta-1}\alpha
 =\binom {\alpha+\beta-1}{\alpha-1} - \binom{\alpha+\beta-1}\alpha.
\end{equation}
This matches the formula $\binom{n}k-\binom{n}{k+1}$ we found for the number
of positive walks of length $n$ ending at~$2k-n\geq0$, with $n=\alpha+\beta-1$
and $k=\alpha-1$.

Clearly that formula is directly related to the arguments of both Bertrand and
Andr{\'e}. Passing from it to the formula for all positive walks just requires
summing over all $k\geq\frac{n}2$, which gives a telescoping sum adding up
to~$\smash{\binom{n}{\lceil{n}/2\rceil}}$. (The same argument shows that there
are $\binom{n}k$ positive walks ending at a number that is \emph{at
  least}~$2k-n$.) It may therefore be expected that there is a close relation
between bijections proving respectively Andr{\'e}'s claim of parity between of
bad and ugly cases, and the enumerative statement of our theorem~1.

Consider our basic raising operation, which as we recall reverses the first
step of a walk that reaches the global minimum of that walk. Within the set of
all walks ending at non-negative numbers, it is defined on the subset of
walks that do descend at least once to a negative number, and its image is the
set of walks that end at a value~$\geq2$. Note that these sets are the
complements of respectively the set of positive walks and the set of walks
that end at $0$~or~$1$, so that we can interpret the fact that iteration of
this raising operation defines a bijection between those complements (in the
opposite direction) as an instance of the ``complementation principle'' evoked
in the introduction.

Upon restriction to the walks that end at a given value~$i\geq0$, the domain
of the raising operation is still limited by the requirement that the walk
descend below~$0$, but the image is the full set of walks ending at~$i+2$. We
can match these two sets with set of ballot sequences for the ugly and bad
cases, provided we omit from those sequences the initial vote (whose outcome
is fixed in either case). Thus we get walks starting at $1$ respectively
at~$-1$, of length $n=\alpha+\beta-1$, and the ending at $\alpha-\beta$; this
point is ahead of the starting point by the margin of $A$ over~$B$ among the
votes remaining after the first one, namely $i=(\alpha-1)-\beta$ in the ugly
cases and $\alpha-(\beta-1)=i+2$ in the bad cases. Since the ugly cases now
start with a one vote lead for~$A$, the condition of subsequently losing this
strict lead at least once matches the condition of descending below the
starting point that characterises the domain of the raising operation. Thus we
now interpret the raising operation as one that fixes the end point, and
shifts back the starting point by $2$~units, as illustrated in figure~3.
\begin{figure}[t]
$$\vcsanepsf{\epsfysize=52mm}{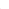}$$
\caption{The raising operator as a starting-point shifting transformation}
\end{figure}
So the raising operation does indeed provide a bijection that proves the claim
in Andr{\'e}'s argument.

This does not of course imply that the raising operation defines the same
bijection that Andr{\'e} did; it does not. We shall give the details of the
latter bijection below, as well as the descriptions of other bijections that
link the sets of bad and ugly ballot sequences, and therefore could be used
instead to complete Andr{\'e}'s argument. But we first want to observe that
\emph{any} such bijection can be iterated to obtain a bijection between almost
recurrent and positive paths: this is not a particularity of the raising
operation.

\begin{prop}
Let a bijection~$f$ be given that maps ugly ballot sequences to bad ones,
without changing the number of votes for either candidate. One can then
bijectively map the set of non-bad ballot sequences in which $A$ has a final
margin of $1$~or~$2$ votes over~$B$ to arbitrary good ballot sequences of the
same length, by iterating the following step as long as the ballot sequence is
ugly: apply~$f$ to change it into a bad sequence, then replace the initial
vote for~$B$ by one for~$A$.
\end{prop}

\emph{Proof.} Termination is ensured since $A$ gains one vote at each
iteration, and bijectivity of the resulting correspondence follows immediately
from that of~$f$. \QED

One may view the proposition as a mere instance of the complementation
principle, within the set of all ballot sequences starting with a vote for~$A$
and resulting in a victory for~$A$. In fact, this shows that the hypothesis of
not changing the number of votes for either candidate is superfluous, as long
as $A$ still wins after applying~$f$; but without it one would not have
such a good upper bound for the number of iterations required as one has with
the hypothesis, namely the initial number of votes for~$B$. From the bijection
obtained by this iteration, an alternative proof for the enumerative statement
of our theorem~1 can be deduced by removing the initial vote (for~$A$) from
each sequence.

Bijections~$f$ as in the proposition are often naturally thought of as
acting on lattice paths, going from the origin by downwards diagonal steps
to a common end point~$(n,i)$, with $i\geq0$; in this point of view the change
of the initial vote (back) into one for~$A$ done at the end of each iteration
shifts the whole remaining path to the right. It might in some cases be more
transparent to present the iteration instead as involving a move of the origin
to the left at each step.

It is clear that like ugly ones, bad ballot sequences must reach a point where
the votes for $A$ and~$B$ are in balance. This circumstance provides an
occasion, when seeking a bijection~$f$ as in the proposition, to focus on the
sequences up to this point of equality, and leave their remainders unchanged.
Moreover the most obvious way to see that there are as many sequences of a
given length ending in equality, but with $A$ ahead of~$B$ at all intermediate
times, as there are with $B$ similarly ahead of~$A$, is to simply interchange
the roles of $A$ and~$B$. In terms of lattice paths this amounts to reflecting
their part up to the first encounter of the vertical line through the origin,
which is the symmetry axis with respect to $A$ and~$B$, in that axis. Using
the resulting bijection is known as the ``reflection method'' for solving the
ballot problem. Surprisingly, and in spite of the fact that it is often
attributed to him, the reflection method is \emph{not} what Andr{\'e} used to
complete his proof either (see~\cite{renault} for details about the intriguing
history of this problem). So let us finally describe the bijection he did
define.

To transform an ugly sequence into a bad one, one may split the former
sequence just before the vote where it turns ugly (the vote for~$B$ after
which the initial lead of~$A$ is levelled); the two parts of the sequence are
then concatenated in the opposite order. The result clearly has the same
distribution of votes, and is bad since it starts with a vote for~$B$. That
the operation is bijective follows from the fact that the initial vote for~$A$
in the ugly sequence has become the last of the votes for~$A$ in the bad
sequence to make the margin of $A$ over~$B$ attain its ultimate value; it can
therefore be located in the result, and the cyclic rearrangement of votes
reversed. If one considers the sequences without their initial vote, as is
more practical for the purpose of iteration, then the description changes
slightly: one should remove the first vote that gives~$B$ a lead over~$A$, and
concatenate the remaining parts of the sequence in opposite order with a vote
for~$A$ inserted in between. Curiously the description Andr{\'e} originally
gave includes the initial vote for~$A$ in the ugly sequence but maps to a bad
sequence deprived of its initial vote for~$B$, so his recipe is: remove the
offending vote, and combine the remaining sequences in the opposite order.

It appears there is ample choice for a bijection proving Andr{\'e}'s claim. In
fact, it is not hard to find more bijections, as variations of the reflection
method. In the lattice path view, one is not obliged to take the first visit
of the axis as the point up to which the reflection is applied (unless it is
the unique visit): one could equally well choose say the second visit (if
possible), or the last one; indeed one could fix any rule that depends only on
properties of the path that are unchanged by the reflection. Alternatively,
reflection in the axis is not the only way to map Dyck paths bijectively to
those staying at the opposite side of the axis: central symmetry with respect
to the midpoint on the axis works as well; this corresponds to reversing the
order of counting votes rather than changing individual votes. In this case,
unlike for the reflection method, one must insist on transforming only the
part up to the \emph{first} encounter of the axis, to ensure that then
result corresponds to a sequence starting with a vote for~$B$.

By our proposition, each of these bijections gives rise to a bijection between
(almost) recurrent walks and positive walks of the same length. We conclude
this note by comparing the various possibilities with respect to the
transparency of the resulting bijection. As indicator we consider
understanding the number of iterations required to turn a given walk into a
positive walk. It is clear that this number will equal half the value of the
ending point of that positive walk (rounded down in the odd-length case), and
so it will in all cases define a statistic with values in
$\{0,1,\ldots,\lfloor\frac n2\rfloor\}$ on the set of (almost) recurrent paths
of length~$n$. Moreover this statistic will have the same distribution as the
``half the ending value rounded down'' statistic has on positive walks of
length~$n$, namely $\binom{n}k-\binom{n}{k-1}$ instances for the
value~$k-\lceil\frac{n}2\rceil\geq0$. For the bijection of theorem~1,
iterating the raising operation, we have seen that this statistic is just the
depth of the walk (minus the value of its global minimum), which is easily
read off.

In contrast, finding the number of iterations required for a given walk when
using the bijection given by the reflection method is not at all easy; as far
as we can see it requires essentially simulating the entire iteration process.
The reason for this difficulty is that the part of the walk affected by each
successive reflection might be either smaller or larger than for the previous
reflection, and as a result, computing where a given path segment will end up
after a certain number of iterations becomes a somewhat messy affair.

This difficulty is present also for the variations of the reflection method we
indicated, with one exception: if we systematically reflect the \emph{largest}
possible part of the path, namely up to the last encounter of the axis, then
successive iterations are guaranteed to affect ever smaller parts of the path.
We can then deduce the following computation of the number of iterations
required for this case: trace the walk (as before taken to not include a step
for the initial vote, and to start at~$0$) backwards in time, seeking first
the last visit to~$-1$ (if any), then the preceding visit to~$1$, then the
preceding visit to~$-1$, and so on alternatingly until no visit of the
indicated kind remains; the requested number is the number of visits found.
For what it's worth, we formulate the fact that this statistic has the
mentioned distribution; this is actually not too hard to see directly.

\begin{corol}
Given $0\leq{d}\leq\lfloor{n/2}\rfloor$, there are
$\binom{n}{\lceil\frac{n}2\rceil+d}$ walks on~$\Z$ of length~$n$ starting
at~$0$ and ending at $n\bmod2\in\{0,1\}$ that make at least $d$ alternating
visits to $1$ and $-1$, the last of which is a visit to~$-1$; the number of
such walks for which $d$ is the length of the longest such sequence of visits
is $\binom{n}{\lceil\frac{n}2\rceil+d}-\binom{n}{\lceil\frac{n}2\rceil+d+1}$.
\end{corol}

For Andr\'e's original method, we find that the number of iterations required
is the same as when using the raising operation, namely the depth of the walk.
Indeed each iteration decreases the depth by~$1$, since the part of the
sequence moved to the end corresponds to a Dyck path, and the following
down-step is removed; coming after the absolute minimum is reached, the
additional up-step that is inserted does not affect the depth. In fact the
downs-steps removed in successive iterations are the same ones as those that
iteration of the raising operation would change into up-steps, although the
latter would proceed in the opposite (rear to front) order. Indeed it is easy
to describe the final walk obtained by iteration using Andr\'e's bijection in
terms of the one of theorem~1: from the latter walk remove the last of the
(up-)steps that were obtained by reversing down-steps, and combine the
remaining parts of the walk in the opposite order, separated by an up-step.
This close relation is remarkable, as it seems unlikely that this kind of
iteration was of any concern to Andr\'e.

In conclusion, although one can mechanically transform bijective solutions to
the ballot theorem into alternative bijections for our theorem~1, this turns
the more obvious solutions (notably the reflection method) into rather opaque
bijections. In contrast, the raising operation as well as Andr\'e's almost
forgotten original method, although less obvious in relation to the ballot
problem, lead to quite transparent bijections after iteration.

\end{document}